\input amstex
\documentstyle{amsppt}
\overfullrule=0pt
\topmatter
\title
Mathematics and Economics of
Leonid Kantorovich
\endtitle
\author  S. S. Kutateladze \endauthor
\date    January 2, 2011\enddate
\address
Sobolev Institute of Mathematics\newline
\indent Novosibirsk 630090, Russia
\endaddress
\email   sskut\@math.nsc.ru
\endemail
\keywords
linear programming, functional analysis, applied mathematics
\endkeywords
\abstract
This is a short overview of the  contribution
of Leonid Kantorovich into the formation of the modern
outlook on the interaction between mathematics and economics.
\endabstract
\endtopmatter

\def\st{\operatorname{st}}
\def\ltd{\operatorname{ltd}}

\def\bnorml{\mathopen{\kern1pt \vrule height6.5pt depth1.5pt width1pt\kern1.5pt}}
\def\bnormr{\mathclose{\kern1pt \vrule height6.5pt depth1.5pt width1pt\kern1pt}}

\document

\begingroup

Leonid  Vital'evich Kantorovich was a renowned mathematician and economist,
 a prodigy and a Nobel prize winner.
These extraordinary circumstances deserve some
attention in their own right. But they can lead hardly to any useful
conclusions for the general audience in view of their extremely low probability.
This is not so with the creative legacy of a human, for what is done
for the others remains unless it is forgotten, ruined, or libeled. 
Recollecting the contributions of  persons to culture, we
preserve their spiritual worlds for the future.

\head
Kantorovich's Path
\endhead

Kantorovich was born in the~family of a~venereologist at St.~Petersburg on January 19, 1912
(January~6, according to the
old Russian style).
In 1926, just  at the age of 14,  he entered St. Petersburg
(then Leningrad) State University (SPSU).  His supervisor was
G.~M.~Fikhtengolts.

After graduation from SPSU in 1930, Kantorovich started
teaching, combining it with intensive scientific research. Already in
1932  he became a~full professor at  the Leningrad Institute of
Civil Engineering and an assistant professor at SPSU.
From 1934 Kantorovich
was a~full professor at his {\it alma mater}.

The main achievements in mathematics belong to the
``Leningrad'' period of Kantorovich's life. In the
1930s he published more papers in pure mathematics
whereas his 1940s are devoted to computational mathematics
in which he was soon appreciated as a~leader in this country.

The letter of Academician  N.~N. Luzin, written on April 29, 1934,
was  found in the personal archive of Kantorovich a few years ago
during preparation of his selected works for publication (see~[1]).

This letter demonstrates the attitude of Luzin, one of the most eminent
and influential mathematicians of that time, to the brilliance of
the young prodigy. Luzin was the founder and leader of the famous
``Lusitania'' school of Muscovites.
He remarked in his letter:

\medskip
\hangindent 10 pt{\eightpoint\it
... you must know my attitude to you.  I do not know you as a man
$\underline{\text{completely}}$
but I guess a warm and admirable personality.

However, one thing I know for certain: the range of your mental powers
which, so far as I  accustomed myself to guess people, open up
$\underline{\text{limitless possibilities}}$  in science.   I will not utter
the appropriate word---what for?
Talent---this would belittle you. You are entitled to get more...
}
\medskip

\noindent
In 1935 Kantorovich made his major mathematical
discovery---he defined $K$-{\it spaces\/}, i.~e., vector lattices whose every nonempty
order bounded subset had an infimum and supremum.
The Kantorovich spaces have provided the natural framework
for developing the theory of linear inequalities
which was a practically uncharted area of research those days.
The concept of inequality is obviously relevant
to approximate calculations
where we are always interested in various estimates
of the accuracy of results.
Another challenging source of interest in linear inequalities
was the stock of problems of economics.
The language of partial comparison is rather natural
in dealing with
what is reasonable and optimal in human behavior
when means and opportunities are scarce.
Finally,
the concept of linear inequality is inseparable
from the key idea of a convex set.
Functional analysis implies
the existence of nontrivial continuous linear functional
over the space under consideration,
while the presence of a functional of this type
amounts to the existence of nonempty proper open convex subset
of the ambient space.
Moreover,
each convex set is generically the solution set of
an appropriate system of simultaneous linear inequalities.

At the end of the 1940s Kantorovich
formulated and explicated the thesis of interdependence between
functional analysis
and applied mathematics:

\medskip
\hangindent 10 pt{\eightpoint\it
There is now a tradition of viewing functional analysis
as a purely theoretical discipline
far removed from direct applications,
a discipline which cannot deal with practical questions.
This article\footnote"${}^{1}$"{Cp.~[2]; the excerpt is taken
from~[3: Part~2, p.~171].}
is an attempt to break with this tradition,
at least to a certain extent,
and to reveal the relationship between functional analysis
and the questions of applied \hbox{mathematics\dots.}}
\medskip

He distinguished the three techniques:
the Cauchy method of majorants also called {\it domination},
the method  of finite-dimensional approximations,
and the Lagrange method for the new optimization problems
motivated by economics.

Kantorovich based his study of the Banach space versions
of the Newton method
on domination in general ordered vector spaces.

Approximation of infinite-dimensional spaces
and operators by their finite-dimen\-sional analogs, which
is {\it discretization},
must be considered alongside the marvelous universal
understanding of computational mathematics as the science of
finite approximations to general (not necessarily metrizable)
compacta.\footnote"${}^2$"{This revolutionary definition  was given in
the joint talk  by S.~L.~Sobolev, L.~A.~Lyusternik,
and L.~V. Kantorovich at the Third All-Union Mathematical Congress
in~1956;  cp. [4, pp.~443--444].}

The novelty of the extremal problems arising in social sciences
is connected with the presence
of multidimensional contradictory utility functions.
This raises the major problem of agreeing  conflicting aims.
The corresponding techniques may be viewed as
an instance of {\it scalarization\/}  of vector-valued targets.

From the end of the 1930s the research of Kantorovich acquired new traits
in his audacious breakthrough to economics.
Kantorovich's booklet
{\it Mathematical Methods in the Organization and Planning of Production\/}
which appeared in 1939 is a material evidence
of the birth of linear programming.
Linear programming is a technique
of maximizing a~linear functional
over the positive solutions of a system of linear inequalities.
It is no wonder that the discovery of linear programming
was immediate after the foundation
of the theory of Kantorovich spaces.

The economic  works of
Kantorovich were hardly visible at the surface of the scientific
information flow in the 1940s. However, the problems of economics
prevailed in his creative studies. During the Second World War  he
completed the first version of his book {\it The Best Use of Economic
Resources} which led to the Nobel Prize awarded to him and
Tjalling C.~Koopmans in~1975.

The Council of Ministers of the USSR issued a top secret
Directive No. 1990--774ss/op\footnote"${}^3$"{The letters ``ss'' abbreviate
the Russian for ``top secret,'' while the letters ``op'' abbreviate
the Russian for ``special folder.''}
in 1948 which  ordered ``to organize in the span of two weeks
a group for computations with the staff up to 15 employees
in the Leningrad Division of the Mathematical Institute of the Academy of Sciences of the USSR
and to appoint Professor Kantorovich the head of the group.''
That was how Kantorovich was enlisted in the squad of
participants of the project of producing  nuclear
weapons in the USSR.\footnote"${}^4$"{This was the Soviet project
 ``Enormous,''  transliterated in Russian like ``\'Enormoz.''
 The code name was used in the  operative correspondence
of the intelligence services of the USSR.}

In 1957 Kantorovich accepted the invitation to join
the newly founded  Siberian Division of the Academy of Sciences of
the USSR.
He  moved to Novosibirsk and soon became a corresponding member
of the Department of Economics in the first elections to the
Siberian Division.
Since then  his major publications were
devoted to economics with the exception of
the celebrated course of functional analysis,
``Kantorovich and Akilov'' in the students' jargon.

It is impossible not to mention one brilliant twist of mind of
Kantorovich and his students in suggesting a scientific approach to
taxicab metered rates.
The people of the elder generation in this country remember that in the 1960s the taxicab meter
rates  were modernized radically: there appeared a price for
taking a taxicab which was combining with a less  per kilometer cost.
This led immediately to raising efficiency of
taxi parks as well as profitability of short taxicab drives.
This economic measure was a result of a mathematical modeling
of taxi park efficiency which was accomplished
by Kantorovich with a group of young mathematicians and
published in the rather prestigious mathematical journal
{\it Russian Mathematical Surveys}.

The 1960s became the decade of his recognition.
In 1964 he was elected a full member of the Department of Mathematics
of the Academy of Sciences of the USSR,
and in 1965 he was awarded the Lenin Prize. In these years he
vigorously propounded and maintained his views of interplay between
mathematics and economics  and exerted great efforts to instill the ideas and methods
of  modern science into the top economic management of the Soviet Union,
which was almost in vain.

At the beginning of the 1970s Kantorovich left Novosibirsk for
Moscow where he was deeply engaged in economic analysis, not
ceasing his efforts to influence  the everyday economic practice and
decision making  in the national economy. His activities were mainly
waste of time and stamina in view of the misunderstanding and
hindrance of the governing retrogradists of this country.
Cancer terminated his path in science on April~7, 1986.
He was buried at Novodevichy Cemetery in Moscow.

\head
Contribution to Science
\endhead

The scientific legacy of Kantorovich is immense. His research in
the areas of functional analysis, computational mathematics,
optimization, and descriptive set theory has had a dramatic impact
on the foundation and progress of these disciplines. Kantorovich
deserved his status of one of the father founders of the modern
economic-mathematical methods. Linear programming, his most
popular and celebrated discovery, has changed the image of economics.

Kantorovich wrote more than 300 articles. When we discussed
with him the first edition of an annotated bibliography of his
publications in the early 1980s, he  suggested to combine them in
the nine sections: descriptive function theory and set theory,
constructive function theory,
approximate methods of analysis,
functional analysis,
functional analysis and applied mathematics,
linear programming,
hardware and software,
optimal planning and optimal prices,
and the economic problems of a~planned economy.

Discussing the mathematical papers of Kantorovich, we must especially
mention the three articles [2, 5, 6] in  {\it Russian Mathematical Surveys}.
The first of them had acquired the title that is still impressive in view of
its scale, all the more if compared with the age of the author.
This article appeared in the formula of the Stalin Prize
of 100,000 Rubles which was awarded to Kantorovich in~1948.
The ideas of this brilliant masterpiece laid grounds for
the classical textbook by Kantorovich and Akilov which was the deskbook of many
scientists of theoretical and applied inclination.

The impressive diversity of these areas of research
rests upon not only  the traits of Kantorovich but also
his methodological views.
He always emphasized the innate
integrity of his scientific research as well as mutual penetration and
synthesis of the methods and techniques he used in solving the
most diverse theoretic and applied problems of mathematics and
economics.

The characteristic feature of the contribution of Kantorovich is
his orientation to the most
topical and difficult problems  of mathematics and economics of his epoch.

\head
Functional Analysis and Applied Mathematics
\endhead

The creative style of Kantorovich rested on the principle of unity
of theoretical and applied studies.
This principle led him to the first-rate achievements at the frontiers
between functional analysis and applied mathematics.
Kantorovich's technique consisted in developing and applying
the methods of domination, approximation, and scalarization.

Let $X$ and $Y$ be real vector spaces lattice-normed
with Kantorovich spaces $E$ and $F$. In other words,
given are some lattice-norms ${\bnorml\cdot\bnormr}_{X}$
and ${\bnorml\cdot\bnormr}_{Y}$. Assume further that $T$
is a linear operator from~$X$ to~$Y$ and $S$ is a~positive
operator from  $E$ into~$F$ satisfying

$$
\CD
X    @>T>>    Y\\
@V{{\bnorml\cdot\bnormr}_{X}}VV       @V{\bnorml\cdot\bnormr}VV\\
E  @>S>>    F
\endCD
$$

Moreover, in case

$$
{\bnorml Tx \bnormr}_{Y} \leq S{\bnorml x\bnormr}_{X}\quad (x\in X),
$$

\noindent
we call $S$  the~{\it dominant\/} or {\it majorant\/} of~$T$.
If the set of all dominants of~ $T$ has the least element, then the latter
is called the {\it abstract norm\/}  or {\it least dominant\/} of $T$ and denoted by
 $\bnorml T\bnormr $.
Hence, the least dominant  $\bnorml T\bnormr $~is the least
positive operator from~$E$ to~ $F$ such that
 $$
 \bnorml Tx\bnormr \leq \bnorml T\bnormr \big(\bnorml  x \bnormr \big)\quad (x\in X).
 $$
\noindent
Kantorovich wrote about this matter in~[7] as
follows:

\medskip
\hangindent 10 pt{\eightpoint\it
The abstract norm enables us to estimate an element much
sharper that a single number, a real norm. We can thus
acquire more precise (and more broad) boundaries of the
range of application of successive approximations.
For instance, as a norm of a continuous function
we can take the set of the suprema of its modulus
in a few partial intervals\dots.
This allows us to estimate the convergence domain of
successive approximations for integral equations.
In the case of an infinite system of equations we know that
each solution  is as a sequence and we can take as the norm of
a sequence not only a sole number but also finitely many numbers;
for instance, the absolute values of the first entries and the estimation
of the remainder:
$$
 \bnorml(\xi_1,\xi_2,\dots)\bnormr=
(|\xi_1|,|\xi_2|,\dots,|\xi_{N-1}|, \sup\limits_{k\geq N}|\xi_k|)\in\Bbb R^N.
$$
}

\hangindent 10 pt{\eightpoint\it
 This enables us to specify the conditions of
applicability of successive approximations
for infinite simultaneous equations. Also, this
approach  allows us to obtain approximate (surplus or deficient)
solutions of the problems under consideration with
simultaneous error estimation.
I believe that the use of
members of semiordered linear spaces instead of reals in various
estimations can lead to essential improvement of the latter.}

\medskip
\noindent
The most general domination underlaid the classical studies of Kantorovich on the Newton method
which brought him international fame.

These days the development of domination proceeds within the frameworks
of Boolean valued analysis (cp.~[8]).
The modern technique of mathematical modeling opened an opportunity
to demonstrate that the principal properties of lattice normed spaces
represent the Boolean valued interpretations of the relevant properties
of classical normed spaces. The most important interrelations here are
as follows: Each Banach space inside a Boolean valued model becomes
a universally complete Banach--Kantorovich space in result of
the external deciphering of constituents. Moreover,
each lattice normed space  may be realized as a dense subspace
of some Banach space in an appropriate Boolean valued model.
Finally, a Banach space~$X$ results from some
Banach space inside  a Boolean valued model
by a special machinery of bounded descent
if and only if $X$ admits a complete Boolean
algebra of norm-one projections which enjoys the
cyclicity property. The latter amounts to the
fact that $X$ is a~Banach--Kantorovich space and  $X$
is furnished with a mixed norm.\footnote"${}^5$"{The modern
theory of dominated operators is thoroughly set forth in
the book~[9] by~A.~G.Kusraev.}

Summarizing his research into the general theory of
approximation methods, Kantorovich wrote:\footnote"${}^6$"{Cp.~[10].}

\medskip
\hangindent 10 pt{\eightpoint\it
There are many distinct methods for various classes of problems
and equations, and constructing and studying them in each particular case presents
considerable difficulties.Therefore, the idea arose of evolving a general
theory that would make it possible to construct and study them with a single source.
This theory was based on the idea of the connection between the given space, in
which the equation to be studied is specified, and a more simple one into which
the initial space is mapped. On the basis of studying the ``approximate
equation'' in the simpler space the possibility of constructing and studying approximate
methods in the initial space was discovered\dots.}

\hangindent 10 pt{\eightpoint\it
It seems to me that the main idea of this theory is of a general character
and reflects the general gnoseological principle for studying
complex systems. It was, of course, used earlier, and it is also used
in systems analysis, but it does not have a rigorous mathematical apparatus.
The principle consists simply  in  the fact that
to a given large complex system in some space a simpler, smaller dimensional
model in this or a~simpler space is associated by means of one-to-one or
one-to-many  correspondence. The study of this simplified model
turns out, naturally, to be simpler and more practicable.
This method, of course, presents definite requirements on the quality
of the approximating system.}
\medskip
\noindent
The classical scheme of discretization
as suggested by Kantorovich for
the analysis of the equation
$
Tx=y,
$
with $T:X\to Y$ a~bounded linear operator
between some Banach spaces $X$ and~$Y$, consists in choosing
finite-dimensional approximating subspaces
$X_N$ and~$Y_N$ and the corresponding embeddings $\imath_N$ and~$\jmath_N$:

$$
\CD
X    @>T>>    Y\\
@AA{\imath_N}A       @A{\jmath_N}AA\\
X_N  @>T_N>>    Y_N
\endCD
$$

\noindent
In this event, the equation
$
 T_N x_N=y_N
$
is viewed as a finite-dimensional approximation
to the original problem.

Boolean valued analysis enables us to
expand the range of applicability of Banach--Kantorovich spaces
and more general modules for studying extensional equations.

Many promising possibilities are open by the new method
of {\it hyperapproximation\/} which rests on the ideas of
infinitesimal analysis.  The classical discretization
approximates  an infinite-dimensional space with the aid of
finite-dimensional subspaces.
Arguing within nonstandard set theory
we may approximate an~infinite-dimensional vector space
with external finite-dimensional spaces.
Undoubtedly, the dimensions of  these hyperapproximations
are given as actually infinite numbers.

The tentative scheme of hyperapproximation
is reflected by the following diagram:

$$
\CD
E    @>T>>    F\\
@V{\varphi_E}VV      @V{\varphi_F}VV\\
E^{\#}  @>T^{\#}>>    F^{\#}
\endCD
$$

\noindent
Here $E$ and $F$  are normed vector space
over the same scalars; $T$ is a~bounded linear operator
from $E$ to~$F$; and ${}^{\scriptscriptstyle\#}$ symbolizes
the taking of the relevant nonstandard hull.

Let $E$~be an~internal vector space over
${}^\ast\Bbb F$, where $\Bbb F$~is the
{\it basic field\/} of scalars; i.~e., the {\it reals}
$\Bbb R$
or {\it complexes}~$\Bbb C$, while  ${}^\ast$ is the symbol of the
Robinsonian standardization.
Hence, we are given the two internal operations
 $+:E\times E\to E$ and $\cdot:{}^\ast\Bbb F\times E\to E$
satisfying the usual axioms of a vector space.
Since $\Bbb F\subset{}^\ast \Bbb F$,  the internal vector space
$E$ is a vector space over $\Bbb F$ as well. In other words,
$E$ is an external vector space which is not a normed nor Hilbert space
externally even if $E$ is endowed with either structure as an internal space.
However, with each normed or pre-Hilbert space we can associate
some external Banach or Hilbert space.

Let $(E,\|\cdot\|)$~be an internal normed space over
~${}^\ast\Bbb F$. As usual, ~$x\in E$
is a~{\it limited\/} element provided that
$\|x\|$~is a limited real
(whose modulus  has a standard upper bound by definition).
If $\|x\|$ is an infinitesimal (=infinitely small real) then  $x$
is also referred to as an {\it infinitesimal}.
Denote by $\ltd(E)$  and  $\mu (E)$  the external sets of limited elements
and infinitesimals of~$E$. The set $\mu (E)$ is the~ {\it monad\/}
of the origin in~$E$.
Clearly, $\ltd(E)$~is an external vector space over~$\Bbb F$, and
 $\mu(E)$~is a subspace of $\ltd(E)$.  Denote the factor-space
 $\ltd(E)/\mu (E)$ by~$E^{\scriptscriptstyle\#}$.
The space  $E^{\scriptscriptstyle\#}$ is endowed with the
natural norm by the formula
$
  \|\varphi x\|:=\|x^{\scriptscriptstyle\#}\|:=\st (\|x\|)\in\Bbb F\quad ({x\in\ltd(E)}).
$
 Here
 $\varphi:=\varphi_E:=(\cdot)^{\scriptscriptstyle\#}:\ltd(E)\to E^{\scriptscriptstyle\#}$
 is the canonical homomorphism, and $\st$ stands for the taking of
the standard part of a limited real.
In this event $(E^{\scriptscriptstyle\#}, {\|\cdot\|})$ becomes an
external normed space that is called the
 {\it nonstandard hull\/} of~$E$.
If $(E,\|\cdot\|)$ is a standard space then
the nonstandard hull of~$E$ is by definition
 the space $({}^\ast\!E)^{\scriptscriptstyle\#}$ corresponding to the
Robinsonian standardization ${}^\ast\!E$.

 If $x\in E$ then $\varphi ({}^\ast x)=({}^\ast x)^{\scriptscriptstyle\#}$
belongs to $({}^\ast\!E)^{\scriptscriptstyle\#}$. Moreover,
 $\|x\|=\|({}^\ast x)^{\scriptscriptstyle\#}\|$.
Therefore, the mapping
 $x\mapsto ({}^\ast x)^{\scriptscriptstyle\#}$ is an isometric embedding
of~$E$ in~$({}^\ast\!E)^{\scriptscriptstyle\#}$.
 It is customary to presume that $E\subset ({}^\ast\!E)^{\scriptscriptstyle\#}$.

 Suppose now that
 $E$ and $F$ are internal normed spaces and
 $T:E \to F$ is an internal bounded linear operator.
The set of reals
 $$
  c(T):= \{C\in {}^\ast{\Bbb R}:\,(\forall x\in E)
 \| T x \| \leq C\|x \| \}
 $$
 is internal and bounded. Recall that
 $\|T\|:= \inf c(T)$.

 If the norm
 $\|T \|$ of~$T$ is limited then the classical normative inequality
 $\| T x \| \leq \| T \|\,\|x\|$ valid for all $x \in E$,
 implies that
 $T({\ltd(E)})\subset \ltd(F)$
and
 $T({\mu (E)})\subset \mu (F)$.
 Consequently, we may soundly define the
descent of~$T$ to the factor space $E^{\scriptscriptstyle\#}$
as the external operator
 $T^{\scriptscriptstyle\#}:E^{\scriptscriptstyle\#}\to
 F^{\scriptscriptstyle\#}$,
 acting by the rule
 $
  T^{\scriptscriptstyle\#} \varphi_E x :=\varphi_F T x\quad (x \in E).
 $
 The operator
 $T^{\scriptscriptstyle\#}$
 is linear  (with respect to the members of~${\Bbb F}$)
 and bounded; moreover, $\|T^{\scriptscriptstyle\#}\| =\st(\|T\|)$.
The operator
 $T^{\scriptscriptstyle\#}$ is called the
 {\it nonstandard hull\/}
of~$T$.
It is worth noting that
 $E^{\scriptscriptstyle\#}$ is automatically
a Banach space for each internal (possible, incomplete)
normed space~$E$.
If the internal dimension of an internal normed
space~$E$ is finite then $E$ is referred to as
a~{\it hyperfinite-dimensional\/} space. To each
normed vector space~$E$
there is a hyperfinite-dimensional subspace
$F\subset {}^\ast\!E$ containing all standard members of
the internal space~
${}^\ast\!E$.

Infinitesimal methods  also provide new
schemes for hyperapproximation of general compact spaces.
As an approximation  to a compact space we may take
an arbitrary internal subset containing all standard elements
of the space under approximation.
Hyperapproximation of the present day stems from
Kantorovich's ideas of discretization.

It was in the 1930s that Kantorovich  engrossed in
the practical problems of decision making.
Inspired by the ideas of functional analysis and order,
Kantorovich attacked these problems in the spirit
of searching for an optimal solution.

Kantorovich observed as far back as in~1948
as follows:\footnote"${}^7$"{Cp.~[7].}:

\medskip
\hangindent 10 pt{\eightpoint\it
Many mathematical and practical problems
lead to the necessity of  finding  ``special'' extrema.
On the one hand, those are boundary extrema when some extremal value
is attained at the boundary of the domain of definition of an argument.
On the other hand, this is the case when the functional to be optimized
is not differentiable. Many problems of these sorts are encountered
in mathematics and its applications, whereas the general methods
turn out ineffective in regard to the problems.
}
\medskip
\noindent
Kantorovich was among the first scientists that formulated
optimality conditions in rather general extremal problems.
We view as classical his approach to the theory of optimal transport
whose center is occupied by the Monge--Kantorovich problem; cp.~[11].

Another particularity of the extremal problems stemming from praxis
consists in the presence of numerous conflicting ends and interests
which are to be harmonized.
In fact, we  encounter the instances of multicriteria  optimization
whose characteristic feature is a vector-valued target.
Seeking for an optimal solution in these circumstances, we
must take into account various contradictory preferences which combine
into  a sole compound aim.  Further more, it is impossible as a rule
to distinguish some particular  scalar target and ignore the rest of
the targets without distorting the original statement of the problem
under study.

The specific difficulties of practical problems
and the necessity of reducing them to numerical calculations
let Kantorovich to pondering over the nature of the reals.
He viewed the members of his $K$-spaces as generalized numbers,
developing the ideas that are now collected around the concept of
scalarization.

In the most general sense, scalarization is reduction to numbers.
Since number is a measure of quantity; therefore,
the idea of scalarization is of importance to mathematics in general.
Kantorovich's studies on scalarization
were primarily connectes with the problems of economics he was
interested in from the very first days of his creative path in science.

\head
Mathematics and Economics
\endhead

Mathematics studies the forms of reasoning. The subject of
economics is the circumstances of human behavior.
Mathematics is abstract and substantive, and the professional
decision of mathematicians do not interfere with  the
life routine of individuals. Economics is concrete
and declarative, and the practical exercises of economists
change the life of individuals substantially.
The aim of mathematics consists in  impeccable truths
and methods for acquiring them.
The aim of economics
is the well-being of an individual and
the way of achieving it.
Mathematics never intervenes into the private life of an individual.
Economics touches his purse and bag.
Immense is the list of striking differences between mathematics
and economics.

Mathematical economics is an innovation of the
twentieth century. It is then when the understanding
appeared that the problems of economics need
a completely new mathematical technique.

{\it Homo sapiens\/} has always been and will stay forever
{\it homo economicus}.
Practical economics for everyone as well as their ancestors
is the arena of common sense.
Common sense is a specific ability
of a~ human
to instantaneous moral judgement.
Understanding is higher than common sense
and reveals itself as the adaptability of behavior.
Understanding is not inherited and so it does nor belong to
the inborn traits of a person.
The unique particularity of  humans
is the ability of sharing their understanding, transforming
evaluations into material and ideal
artefacts.

Culture is the treasure-trove of understanding.
The inventory of culture is the essence of
outlook. Common sense is subjective and affine to the
divine revelation of faith that is the force
surpassing the power of external proofs by fact and formal logic.
The verification  of statements with facts and by logic
is a critical  process liberating a human from the errors of
subjectivity.
Science is an unpaved road to objective understanding.
The religious and scientific versions of outlook
differ actually in the methods of codifying  the
artefacts of understanding.

The rise of science as an instrument of understanding
is a long and complicated process.
The birth of ordinal counting
is fixed with the palaeolithic findings hat
separated by hundreds of centuries from the appearance of
a knowing and economic human.
Economic practice precedes
the prehistory of mathematics
that became the science of provable calculations in Ancient Greece
about 2500 years ago.

It was rather recently that the purposeful behavior of humans
under the conditions of
limited resources became the object of science.
The generally accepted date of the birth of economics as a science is
March 9, 1776---the day when there was published
the famous book by Adam Smith
{\it An Inquiry into the Nature and Causes of the Wealth of Nations}.

\head
Consolidation of  Mind
\endhead

Ideas rule the world.
John Maynard Keynes completed this banal statement with a touch
of bitter irony.
He finished his most acclaimed treatise {\it The General Theory of Employment, Interest,
and Money\/} in a rather aphoristic manner:

\medskip
\hangindent 10 pt{\eightpoint\it
 Practical men, who
believe themselves to be quite exempt from any intellectual influences, are usually the slaves of some
defunct economist.}
\medskip

Political ideas aim at power, whereas economic ideas aim at
freedom from any power.
Political economy is inseparable from not only the economic practice but also
the practical policy. The political content of economic teachings implies their
special location within the world science.
Changes in epochs, including their technological achievements  and
political utilities, lead to the universal proliferation of spread of the
emotional attitude to economic theories, which drives economics in the position
unbelievable for the other sciences.
Alongside noble reasons for that, there is one rather cynical:
although the achievements of exact sciences drastically change the life of the mankind,
they never touch the common mentality of  humans as vividly and sharply as any statement about
their purses and limitations of freedom.

Georg Cantor, the creator of set  theory,
remarked as far back as in 1883 that ``the essence of mathematics lies entirely
in its freedom.''
The freedom of mathematics does not reduce to
the absence of exogenic restriction on the objects and methods of research.
The freedom of mathematics reveals itself  mostly in the new
intellectual tools for conquering the ambient universe  which are provided by mathematics
for liberation of humans by  widening the frontiers of their independence.
Mathematization of economics is the unavoidable stage of the
journey of the mankind into the realm of freedom.

The nineteenth century is marked with the first attempts at
applying mathematical methods to economics in the
research by   Antoine Augustin Cournot, Carl Marx,
William Stanley Jevons, L\'eon Walras, and his successor in Lausanne University
Vilfredo Pareto.

John von Neumann and Leonid Kantorovich, mathematicians of the first calibre, addressed the economic problems
in the twentieth century.
The former developed game theory, making it an apparatus for the study of economic behavior.
The latter invented linear programming for decision making in the problems of
best use of limited resources.
These contributions of von Neumann and Kantorovich
occupy an exceptional place in science.
They demonstrated that
the modern mathematics opens up,broad opportunities for economic analysis
of practical problems. Economics has been drifted closer to mathematics.
Still remaining a humanitarian science, it mathematizes rapidly,
demonstrating high self-criticism  and an extraordinary
ability of objective thinking.

The turn in the mentality of the mankind that was effected by
von Neumann and Kantorovich  is not always comprehended to full extent.
There are principal distinctions between the exact and humanitarian
styles of thinking. Humans are prone to reasoning by analogy
and using  incomplete induction, which invokes the illusion of
the universal value of the tricks we are accustomed to.
The differences in  scientific technologies are not distinguished
overtly, which in turn contributes to self-isolation and deterioration of
the vast  sections of science.

The methodological precipice  between
economists and mathematics was well described
by Alfred Marshall, the founder of the Cambridge school of neoclassicals,
``Marshallians.'' He wrote in his {\it magnum opus}~[12]:

\medskip
\hangindent 10pt {\eightpoint\it
The function then of analysis and deduction in economics
is not to forge a few long chains of reasoning, but to forge
rightly many short chains and single connecting
links...}\footnote"${}^8$"{Cp.~[12, Appendix C: The Scope and Method of Economics. \S~3.]}

\medskip
\hangindent 10pt {\eightpoint\it
It is obvious that there is no room in economics
for long trains of deductive reasoning.}\footnote"${}^9$"{Cp.~[12, Appendix~D:
Use of Abstract Reasoning in Economics].}
\medskip

\noindent
In 1906 Marshall formulated his scepticism in regard to
mathematics as follows:

\medskip
\hangindent 10pt {\eightpoint\it
[I had] a growing feeling in the
later years of my work at the subject that a good mathematical theorem
dealing with economic hypotheses was very unlikely to be good
economics: and I went more and more on the rules---}

\hangindent 30pt {\eightpoint\it
(1) Use mathematics
as a shorthand language, rather than an engine of inquiry.}

\hangindent 30pt {\eightpoint\it
(2) Keep to  them till you have done.}

\hangindent 30pt {\eightpoint\it
(3) Translate into English.}

\hangindent 30pt {\eightpoint\it
(4) Then illustrate by examples that are important in real life.}

\hangindent 30pt {\eightpoint\it
(5) Burn the  mathematics.}

\hangindent 30pt {\eightpoint\it
(6) If you can't succeed in (4), burn (3). This last I
did often.}

\hangindent 30pt {\eightpoint\it
I don't mind the mathematics, it's useful and necessary,
but it's too bad the history of economic thought is no longer required
or even offered in many graduate and undergraduate programs. That's a
loss. }\footnote"${}^{10}$"{Cp.~[13, p.~294].}
\medskip

\noindent
Marshall intentionally counterpose the economic
and mathematical ways of thinking, noting that
the numerous short ``combs'' are appropriate in a concrete economic analysis.
Clearly, the image of a~``comb'' has nothing in common with
the upside-down pyramid, the cumulative hierarchy of
the von Neumann universe, the residence of the modern
Zermelo--Fraenkel set theory.
It is from  the times of  Hellas
that the beauty and power of mathematics rest on
the axiomatic method which presumes the derivation of new facts by
however lengthy chains of formal implications.

The conspicuous discrepancy between economists and mathematicians
in mentality has hindered their mutual understanding and cooperation.
Many partitions, invisible  but ubiquitous, were erected in ratiocination,
isolating  the economic community from its mathematical
counterpart and vice versa.

This status quo  with deep roots in history
was always a~challenge to Kantorovich, contradicting
his views of interaction between mathematics and economics.

\head
Linear Programming
\endhead

The principal discovery of Kantorovich
at the junction of mathematics and economics is linear programming which is now
studied by hundreds of thousands of people throughout the world.
The term signifies the colossal area of science which is allotted to
linear optimization models.
In other words, linear programming is the science of the theoretical
and numerical analysis of  the problems in which we seek for an optimal (i.~e.,
 maximum or minimum) value of some system of indices of a process
 whose behavior is described by simultaneous linear inequalities.

The term ``linear programming'' was minted in 1951
by  Koopmans. The most commendable contribution of Koopmans
was the ardent promotion of the methods of linear programming and
the strong defence of Kantorovich's priority in the invention of these methods.

In the USA the independent research into linear optimization models was started only
in 1947 by George B. Dantzig who convincingly described the history of the area
in his classical book [14, p.~22--23] as follows:

\medskip
\hangindent 10pt {\eightpoint\it
The  Russian  mathematician L.~V.  Kantorovich  has  for  a
number  of  years  been  interested in  the  application  of
mathematics   to  programming  problems.  He  published   an
extensive monograph in 1939 entitled {\sl Mathematical Methods in
the Organization and Planning of Production}}...

\hangindent 10pt {\eightpoint\it
Kantorovich  should  be credited with being  the  first  to
recognize that certain important broad classes of production
problems had well-defined mathematical structures which,  he
believed,  were  amenable to practical numerical  evaluation
and could be numerically solved.}

\hangindent 10pt {\eightpoint\it
In  the  first  part of his work Kantorovich  is  concerned
with  what  we  now call the weighted two-index distribution
problems. These were generalized first to include  a  single
linear  side  condition,  then  a  class  of  problems  with
processes     having     several    simultaneous     outputs
(mathematically the latter is equivalent to a general linear
program). He outlined a solution approach based on having on
hand  an  initial feasible solution to the  dual.  (For  the
particular problems studied, the latter did not present  any
difficulty.)  Although the dual variables  were  not  called
``prices,''  the general idea is that the assigned  values  of
these  ``resolving multipliers'' for resources in short supply
can  be  increased  to a point where it  pays  to  shift  to
resources that are in surplus. Kantorovich showed on  simple
examples  how  to make the shifts to surplus  resources.  In
general,  however, how to shift turns out  to  be  a  linear
program  in  itself  for which no computational  method  was
given.  The  report  contains an outstanding  collection  of
potential applications...}

\hangindent 10pt {\eightpoint\it
If  Kantorovich's earlier efforts had been  appreciated  at
the  time  they  were first presented, it is  possible  that
linear  programming would be more advanced  today.  However,
his  early work in this field remained unknown both  in  the
Soviet  Union  and  elsewhere for nearly two  decades  while
linear  programming became a highly developed art.}
\medskip

\noindent
It is worth observing that to an optimal plan of every
linear program there corresponds some optimal prices
or ``objectively determined estimators.'' Kantorovich
invented this  bulky term by tactical reasons in order to
enhance the ``criticism endurability'' of  the concept.

The interdependence of optimal solutions and optimal prices is
the crux of the economic discovery of Kantorovich.

\head
Universal Heuristics
\endhead

The integrity of the outlook of Kantorovich was revealed in all instances of
his versatile research. The ideas of linear programming were tightly interwoven
with his methodological standpoints in the realm of mathematics.
Kantorovich viewed as his main achievement in this area
the distinguishing of
{\it $K$-spaces\/}.\footnote"${}^{11}$"{Kantorovich wrote about ``my spaces'' in his personal memos.}

Kantorovich observed in his first  short paper
of 1935 in {\it Doklady\/} on the newly-born area of
ordered vector spaces:\footnote"${}^{12}$"{Cp.~[3: Part~2, pp.~213--216].}

\medskip
\hangindent 10pt {\eightpoint\it
In this note, I define a~new type of space that I call
a~semiordered linear space. The introduction of such a~space allows us to
study linear operations of one abstract class
(those with values in these spaces) in the same way as linear functionals.}
\medskip

\noindent
This was the first formulation of the most important
methodological position that is now referred to
{\it Kantorovich's heuristic principle}.
It is worth noting that his definition of a semiordered
linear space contains the axiom of Dedekind completeness which
was denoted by $I_6$.  Kantorovich demonstrated the role of
$K$-spaces by widening  the scope of the  Hahn--Banach Theorem.
The heuristic principle turned out applicable to this fundamental
Dominated Extension Theorem; i.e., we may abstract the Hahn--Banach Theorem on
substituting the elements of an arbitrary  $K$-space
for reals  and replacing linear functionals
with operators acting into the space.

Attempts at formalizing Kantorovich's heuristic
principle  started in the middle of the twentieth century
at the initial stages of $K$-space theory and yielded the so-called
identity preservation theorems.
They assert that if some algebraic proposition with
finitely many function variables is satisfied by  the assignment of
all real values then
it remains valid  after replacement of reals with members of an~
arbitrary $K$-space.

Unfortunately,  no satisfactory explanation
was suggested for the internal mechanism behind the phenomenon
of identity preservation. Rather obscure remained
the limits on the heuristic transfer principle. The same applies
to  the general reasons for similarity
and parallelism between the reals  and their analogs
in $K$-space which reveal themselves every now  and then.

 The abstract theory of $K$-spaces, linear programming, and
approximate methods of analysis were particular outputs of Kantorovich's
universal heuristics.
In his last mathematical paper [15] which was written when he had been mortally ill,
Kantorovich remarked:

\medskip
\hangindent 10pt {\eightpoint\it
One aspect of reality was temporarily omitted in the
development of the theory of function spaces.
Of great importance  is the relation of comparison
between practical objects, alongside algebraic and other
relations between them. Simple comparison applicable to every pair
of objects is of a depleted character; for instance,
we may order all items by weight which is of little avail.
That type of ordering is more natural
which is defined or distinguished when this is reasonable
and which is left indefinite otherwise
 ({\sl partial ordering or semiorder}).
For instance,
two sets of goods must undoubtedly be considered as comparable
and one greater than the other if
each item of the former set is quantitatively greater than
its counterpart in the latter. If some part of the goods of
one set is greater  and another part is less than
the corresponding part of the other then
we can avoid prescribing any order
between these sets.  It is with this in mind that
the theory of ordered vector spaces was propounded
and, in particular, the theory of the above-defined  $K$-spaces.
It found various applications not only in
the theoretic problems of analysis but also
in construction of some applied methods, for instance
the theory of majorants in connection
with the study of successive approximations.
At the same time the opportunities it offers are not
revealed fully yet.
The importance  for economics is underestimated of this branch of
functional  analysis.
However, the comparison and correspondence relations
play an extraordinary role in economics and it was
definitely clear even at the cradle of $K$-spaces
that they will find their place in economic analysis and
yield luscious fruits.}

\smallskip
\hangindent 10pt {\eightpoint\it
The theory of $K$-spaces has another important feature:
their elements can be treated as numbers.
In particular, we may use elements of such a space, finite- or
infinite-dimensional, as a norm in construction of analogs of  Banach spaces.
This choice of norms for objects is much more accurate.
Say, a function is normed not by its maximum on the whole
interval but a dozen of numbers, its maxima on parts of this interval}.
\medskip

\noindent
More recent research has corroborated that
the ideas of linear programming
are immanent in the theory of  $K$-spaces. It was demonstrated
that the validity of one of the various statements of the duality principle of
linear programming in an abstract mathematical structure implies with necessity
that the structure under consideration is in fact a~$K$-space.

The Kantorovich heuristic principle is connected with one of the
most brilliant pages of the mathematics of the twentieth century---
the famous problem of the continuum.
Recall that some set $A$ has the cardinality of the continuum
whenever $A$ in equipollent with a segment of the real axis.
The continuum hypothesis is that each subset of the segment is either
countable of has the cardinality of the continuum.
The continuum problem asks whether the continuum hypothesis is true or false.

The continuum hypothesis was first conjectured by Cantor in~1878.
He was convinced that the hypothesis was a theorem and vainly attempted at proving
it during his whole life. In~1900 th Second Congress of Mathematicians took place in
Paris. At the opening session Hilbert delivered his epoch-making talk
``Mathematical Problems.'' He raised  23~problems whose solution was
the  task of the nineteenth century bequeathed to the twentieth century.
The  first on the Hilbert list was open  the continuum problem.
Remaining unsolved for decades, it gave rise
to deep foundational studies.
The efforts of more than a half-century
yielded the solution: we know now that the continuum hypothesis
can neither be proved nor refuted.

The two stages led to  the understanding that the continuum hypothesis is an independent axiom.
G\"odel showed in 1939 that the continuum hypothesis is consistent with
the axioms of set theory\footnote"${}^{13}$"{This was done for ZFC.}, and Cohen demonstrated in
1963 that the negation of the continuum hypothesis does not contradict
the axioms of set theory either.
Both results were established by exhibiting appropriate models; i.~e.,
constructing a universe and interpreting set theory in the universe.
 The G\"odel approach based on ``truncating'
  the von Neumann universe.
G\"odel proved that the constructible sets he distinguished
yield the model that satisfies the continuum hypothesis.
Therefore, the negation of the continuum hypothesis is not provable.
The approach by Cohen was in a sense opposite to that of G\"odel:
it used a controlled enrichment of the von Neumann universe.

Cohen's method of  forcing  was simplified in~1965
on using the tools of Boolean algebra and the new technique
of mathematical modeling which is based on the nonstandard
models of set theory.  The progress of the so-invoked Boolean valued analysis
has demonstrated the fundamental importance of the so-called universally complete
$K$-spaces. Each of these spaces turns out to present
one of the possible noble models of the real axis and so such a space
plays a similar key role in mathematics. The spaces of Kantorovich
implement new models of the reals, this earning their eternal immortality.

Kantorovich heuristics has received brilliant corroboration, this proving
the integrity of science and inevitability of interpenetration of
mathematics and economics.

\head
Memes for  the Future
\endhead

The memes of Kantorovich have been received as witnessed by the curricula and syllabi
of every economics department in any major university throughout the world.
The gadgets of mathematics  and the idea of optimality belong
to the tool-kit of any practicing economist.
The new methods erected an unsurmountable firewall
against the traditionalists that  view economics as a testing polygon
for the technologies like Machiavellianism,
flattery, common sense, or foresight.

Economics as an eternal boon companion of mathematics will avoid merging
into any esoteric part  of~the humanities, or politics, or belles-lettres.
The new generations of mathematicians will treat the puzzling problems
of economics as an inexhaustible source of inspiration and
an attractive arena for applying and refining their impeccably rigorous methods.

Calculation will supersede prophesy.

\endgroup

\enddocument

\end